\documentclass{article}

\usepackage{
 amsmath,
 amsxtra,
 amsthm,
 amssymb,
 etex,
 mathrsfs,
 mathtools,
  stmaryrd
 }
\usepackage[all]{xy}
\usepackage{hyperref}

\newtheorem{theorem}{Theorem}[section]
\newtheorem{lemma}[theorem]{Lemma}

\newtheorem{proposition}[theorem]{Proposition}
\newtheorem{corollary}[theorem]{Corollary}

\newtheorem{lthm}{Theorem} 

\theoremstyle{remark}
\newtheorem{remark}[theorem]{Remark}

\newtheorem{defn}[theorem]{Definition}

\newcommand{\bd}{\begin{defn}}
\newcommand{\ed}{\end{defn}}
\newcommand{\bl}{\begin{lemma}}
\newcommand{\el}{\end{lemma}}
\newcommand{\bp}{\begin{proposition}}
\newcommand{\ep}{\end{proposition}}
\newcommand{\bt}{\begin{theorem}}
\newcommand{\et}{\end{theorem}}
\newcommand{\bc}{\begin{corollary}}
\newcommand{\ec}{\end{corollary}}
\newcommand{\br}{\begin{remark}}
\newcommand{\er}{\end{remark}}
\newcommand{\ba}{\begin{array}}
\newcommand{\ea}{\end{array}}
\newcommand{\bpf}{\begin{proof}}
\newcommand{\epf}{\end{proof}}

\newcommand{\Q}{\mathbb{Q}}
\newcommand{\Qp}{\mathbb{Q}_p}
\newcommand{\Z}{\mathbb{Z}}
\newcommand{\Zp}{\mathbb{Z}_p}
\newcommand{\Op}{\mathcal{O}}
\newcommand{\Ga}{\Gamma}

\newcommand{\la}{\lambda}

\DeclareMathOperator{\Gal}{Gal}
\DeclareMathOperator{\Hom}{Hom}

\DeclareMathOperator{\rank}{rank}

\DeclareMathOperator{\Ext}{Ext}
\newcommand{\Iw}{\mathrm{Iw}}

\newcommand{\lra}{\longrightarrow}
\newcommand{\ord}{\mathrm{ord}}
\newcommand{\ot}{\otimes}
\newcommand{\Tr}{\mathrm{Tr}}

\newcommand{\cyc}{\mathrm{cyc}}
\newcommand{\ps}[1]{\llbracket #1 \rrbracket}

\newcommand{\ilim}{\displaystyle \mathop{\varinjlim}\limits}
\newcommand{\plim}{\displaystyle \mathop{\varprojlim}\limits}

\numberwithin{equation}{section}

\begin{document}
\title{On even $K$-groups over $p$-adic Lie extensions of global function fields}
 \author{
  Meng Fai Lim\footnote{School of Mathematics and Statistics,
Central China Normal University, Wuhan, 430079, P.R.China.
 E-mail: \texttt{limmf@ccnu.edu.cn}} }
\date{}
\maketitle

\begin{abstract} \footnotesize
\noindent Let $p$ be a fixed prime number, and let $F$ be a global function field with characteristic not equal to $p$. In this paper, we shall study the variation properties of the Sylow $p$-subgroups of the even $K$-groups in a $p$-adic Lie extension of $F$. When the $p$-adic Lie extension is assumed to contain the cyclotomic $\Zp$-extension of $F$, we obtain growth estimate of these groups. We also establish a duality between the direct limit and inverse limit of the even $K$-groups.


\end{abstract}

\section{Introduction}

Throughout this article we fix a prime $p$ and an integer $i\geq 2$. Let $F$ be a global function field whose characteristic is not equal to our prime $p$. Once and for all, we fix a separable closure $F^{\mathrm{sep}}$ of $F$. Let $S$ denote a nonempty finite set of primes of $F$, and write $F_S$ for the maximal extension of $F$ contained in $F^{\mathrm{sep}}$ that is unramified outside the set $S$. For every finite extension $L$ of $F$ contained in $F_S$, we write $S(L)$ for the set of primes of $L$ above $F$, and for brevity, we write $\Op_{L,S}$ for the ring of $S(L)$-integers of the field $L$. We also write $G_S(L)= \Gal(F_S/L)$. The goal of this article is to study the even $K$-groups $K_{2i-2}(\Op_{L,S})[p^\infty]$ as $L$ varies in a given $p$-adic Lie extension $F_\infty$ of $F$.

Our approach is via cohomology. Indeed, thanks to the deep work of Rost-Voevodsky \cite{Vo}, the group $K_{2i-2}(\Op_{L,S})[p^\infty]$ can be identified with the continuous cohomology
group $H^2\big(G_S(L), \Zp(i)\big)$. Following standard Iwasawa theoretical procedure, we consider the inverse limit of these cohomology groups, and the resulting object is usually called the Iwasawa cohomology group. More precisely, the $k$-th Iwasawa cohomology group of $\Zp(i)$ over $F_\infty/F$ is defined by
 \[H^k_{\Iw,S}\big(F_\infty/F, \Zp(i)\big):= \plim_L H^k\big(G_S(L), \Zp(i)\big), \]
 where $L$ runs through all finite extensions of $F$ contained in $F_\infty$ and the transition maps are given by corestriction on cohomology. These Iwasawa cohomology groups come naturally equipped with a $\Zp\ps{G}$-module structure, where $G=\Gal(F_\infty/F)$. Furthermore, it can be shown that they are finitely generated over $\Zp\ps{G}$ (for instance, see \cite[Proposition 4.1.3]{LimSh}) and that $H^k_{\Iw,S}\big(F_\infty/F, \Zp(i)\big)=0$ for $k\neq 1,2$. For the first and second Iwasawa cohomology groups, we can even say more and this is the content of the following theorem.

\begin{lthm}[Theorem \ref{H2 mu=0}]
Retain the setting as above. In particular, we let $i\geq 2$. Then the following statements are valid.
\begin{enumerate}
  \item[$(a)$] The module $H^2_{\Iw,S}(F_\infty/F, \Zp(i))$ is torsion over $\Zp\ps{G}$. In the event that $F_\infty$ contains the cyclotomic $\Zp$-extension $F_\cyc$ of $F$, the module $H^2_{\Iw,S}(F_\infty/F, \Zp(i))$ is finitely generated over $\Zp\ps{H}$, where $H= \Gal(F_\infty/F_\cyc)$.
  \item[$(b)$] $H^1_{\Iw,S}(F_\infty/F, \Zp(i))$ is a torsion $\Zp\ps{G}$-module. Furthermore, in the event that $G$ has dimension $\geq 2$, one has
     \[ H^1_{\Iw,S}(F_\infty/F, \Zp(i)) =0. \]
\end{enumerate}
\end{lthm}

We recall that the cyclotomic $\Zp$-extension $F_\cyc$ of $F$ is the unique $\Zp$-extension of $F$ contained in the extension of $F$ obtained by adjoining all $p$-power roots of unity.

The proof of the preceding theorem, which will be given in Section \ref{Iwasawa cohomology groups}, utilizes the fact that the $p$-cohomological dimension of $\Gal(G_S(F)/F_\cyc)$ is equal to one (see \cite[Theorem 10.1.2(c)]{NSW}; also see Lemma \ref{cdpcyc} of this paper). In \cite{DKLL}, an analogous result for the Iwasawa cohomology of the Tate module of an abelian variety is proven by the author and his collaborators. The proof of Theorem \ref{H2 mu=0} will follow that in \cite{DKLL} with some slight modification.

The remainder of the article is concerned with applying the above theorem to the study of even $K$-groups. As a start, the vanishing of the first Iwasawa cohomology group is crucial in the following evaluation of the $G$-Euler characteristics of the second Iwasawa cohomology group.

\begin{lthm}[Theorem \ref{Eulerchar}]
Notation as before. Suppose further that $G=\Gal(F_\infty/F)$ has dimension $\geq 2$ and has no $p$-torsion. Then one has
\[ \chi\big(G,H^2_{\Iw}\big(F_\infty/F, \Zp(i)\big)\big) = \frac{\#K_{2i-2}(\Op_{F,S})[p^\infty]}{\#K_{2i-1}(\Op_{F,S})[p^\infty]}.\]
\end{lthm}
Note that it is a consequence of a theorem of Harder (see \cite[Korollar 3.2.3]{Ha} or \cite[Chap VI, Corollary 6.2]{WeiKbook}) that the group $K_n(\Op_{F,S})$ is finite for $n\geq 2$, and so the above expression on the right makes sense.

The next application of Theorem \ref{H2 mu=0} is concerned with the study of growth formula of even $K$-groups as done in \cite{LimKgrowth}.  However, note that a global function field $F$ with $\mathrm{char}(F)\neq p$ has only one $\Zp$-extension, namely, the cyclotomic one (see \cite[Proposition 10.3.20]{NSW}). Therefore, a $p$-adic Lie extension of dimension $\geq 2$ is necessarily noncommutative. Due to a lack of a precise structural theorem for modules over a noncommutative Iwasawa algebra, this causes an extra complication when dealing with $p$-adic Lie extension with dimension $\geq 2$. Thankfully, in the event that the module is finitely generated over $\Zp\ps{H}$, the author has obtained growth estimates for such modules (see \cite{LimUpper, LimKgrowth}) which in turn are variants of previous results of Perbet \cite{Per} and Lei \cite{Lei} (also see Propositions \ref{ZprtimeZp} and \ref{G/H gen}, and the references given). When $F_\infty$ contains $F_\cyc$, the extra finite generation property of $H^2_{\Iw,S}(F_\infty/F, \Zp(i))$ as seen in Theorem \ref{H2 mu=0} allows to apply the results mentioned to obtain growth estimate. As illustration, we present here the case when $\Gal(F_\infty/F)\cong \Zp^r\rtimes \Zp$.

\begin{lthm}[Theorem \ref{uniform asymptotic arith2}]
Let $i\geq 2$ be given.
  Suppose that $F_\infty$ is a $p$-adic Lie extension of $F$ containing $F^\cyc$ with $G:=\Gal(F_\infty/F_\cyc)\cong \Zp^{d-1}$. Let $F_n$ be the intermediate subfield of $F_\infty/F$ with $\Gal(F_\infty/F_n) = G_n$.
 Then we have
 \[ \ord_p\Big(K_{2i-2}(\Op_{F_n,S})[p^\infty]\Big) = \rank_{\Zp\ps{H}}\Big(H^2_{\Iw,S}\big(F_\infty/F, \Zp(i)\big)\Big) np^{(d-1)n} +  O(p^{(d-1)n}). \]
\end{lthm}

Here $G_n$ denotes the lower $p$-series of $G$ in the sense of \cite{DSMS} (see Section \ref{Growth results section} for the precise definition). For more detailed discussion on growth formula, we refer readers to Section \ref{Growth results section}.

Finally, we come to study the direct limit and inverse limit
\[  \ilim_L K_{2i-2}(\Op_{L,S})[p^\infty] \quad \mbox{and}\quad \plim_L K_{2i-2}(\Op_{L,S})[p^\infty],\]
where the transition maps in the direct limit are induced by the inclusion $\Op_{L,S}\lra \Op_{L',S}$ for $L\subseteq L'$, and the transition maps in the inverse limit are the norm maps (also called the transfer maps). Our result is as follow.

\begin{lthm}[Theorem \ref{Klimitthm}]
Let $i\geq 2$ and let $F_\infty$ be a $p$-adic Lie extension of $F$ which is contained in $F_S$. Suppose that $G=\Gal(F_\infty/F)$ has no $p$-torsion. Then there is an isomorphism
\[  \left(\ilim_L K_{2i-2}(\Op_{L,S})[p^\infty]\right)^\vee \cong \Ext^1_{\Zp\ps{G}}\left(\plim_L K_{2i-2}(\Op_{L,S})[p^\infty], \Zp\ps{G}\right)\]
of $\Zp\ps{G}$-modules.
\end{lthm}

The above result can be thought as an analogue of a classical theorem of Iwasawa \cite{Iw73} which compares the direct limit and inverse of class groups over a $\Zp$-extension of number fields. This result of Iwasawa has been generalized to the context of a $\Zp^d$-extension (see the works of Nekov\'a\v{r} \cite{Ne}, Vauclair \cite{Vau} and, more recently, that of Lai and Tan \cite{LT}). Motivated by these development, the author proved an analogous result for even $K$-groups over a $\Zp^d$-extension of a number field (see \cite{LimKlimit}). A crucial input in the author's proof is the utilization of an exact sequence of Jannsen's (cf. \cite[Theorem 2.1]{Jannsen89} or \cite[Theorem 5.4.13]{NSW}) which yields the following short exact sequence
\[
  0 \lra \Big(\ilim_L K_{2i-2}(\Op_{L,S})[p^\infty]\Big)^\vee \lra  \Ext^1_{\Zp\ps{G}}\left(\plim_L K_{2i-2}(\Op_{L,S})[p^\infty], \Zp\ps{G}\right) \]
  \[  \lra \Big(\ilim_U H_1\big(U,H^2_{\Iw,S}(F_\infty/F, \Zp(i))\big)\ot\Qp/\Zp\Big)^\vee \lra 0
\]
(see Section \ref{Klimitthm section} for details). Therefore, in view of the above exact sequence, it remains to show that the rightmost term vanishes. For the latter, it suffices to show that
$H_1\big(U,H^2_{\Iw,S}(F_\infty/F, \Zp(i))\big)$ is finite. Had $U$ been a commutative group, the finiteness of $H_1\big(U,H^2_{\Iw,S}(F_\infty/F, \Zp(i))\big)$ then follows from the finiteness of $H_0\big(U,H^2_{\Iw,S}(F_\infty/F, \Zp(i))\big)\cong K_{2i-2}(\Op_{L,S})[p^\infty]$ and a standard commutative algebra result (see \cite[P. 56-57]{Se}; also see \cite[Proof of Theorem 2.3]{LimKlimit}). However, in the function field context, every $p$-adic Lie extension of dimension $\geq 2$ is necessarily noncommutative, and so our proof will take a different route. We shall refer readers to Theorem \ref{Klimitthm} for details but we mention here that an important ingredient towards our proof is the finiteness of higher odd $K$-groups. Note that this latter finiteness is not true in general in the number field context.

\subsection*{Acknowledgement}
The author would like to thank Li-Tong Deng, Yukako Kezuka and Yongxiong Li for their interest and comments. In particular, the joint work \cite{DKLL} with them has been a great inspiration to the writing of this paper. The author also likes to thank Chao Qin for helpful comments and discussions relating to the paper. Finally, the author thanks the anonymous referee for the valuable comments and suggestions which helps to improve the clarity of the manuscript.

\section{Iwasawa algebras and modules} \label{Iwasawa modules}

Throughout the paper, $p$ will always denote a fixed prime.
Let $G$ be a compact $p$-adic Lie group. The Iwasawa algebra of $G$ over $\Zp$ is defined by
 \[ \Zp\ps{G} = \plim_U \Zp[G/U], \]
where $U$ runs over the open normal subgroups of $G$ and the inverse
limit is taken with respect to the canonical projection maps.

For now, let us assume that the group $G$ is pro-$p$ and has no $p$-torsion. It is then well-known that $\Zp\ps{G}$ is
a Noetherian Auslander regular ring (cf.\ \cite[Theorem 3.26]{V02}) with no zero divisors (cf.\
\cite{Neu}). As a consequence, it admits a skew field $Q(G)$ which is flat
over $\Zp\ps{G}$ (see \cite[Chapters 6 and 10]{GW} or \cite[Chapter
4, \S 9 and \S 10]{Lam}). This in turn enables us to define the $\Zp\ps{G}$-rank of a finitely generated $\Zp\ps{G}$-module $M$, which is given by
$$ \rank_{\Zp\ps{G}}(M)  = \dim_{Q(G)} (Q(G)\ot_{\Zp\ps{G}}M). $$
The $\Zp\ps{G}$-module $M$ is said to be a
torsion $\Zp\ps{G}$-module if $\rank_{\Zp\ps{G}} (M) = 0$. We also recall a useful equivalent formulation of a torsion $\Zp\ps{G}$-module: $\Hom_{\Zp\ps{G}}(M,\Zp\ps{G})=0$ (cf. \cite[Remark 3.7]{V02}). In the event that the torsion $\Zp\ps{G}$-module $M$ satisfies $\Ext^1_{\Zp\ps{G}}(M,\Zp\ps{G})=0$, we say that $M$ is a pseudo-null $\Zp\ps{G}$-module.
For a finitely generated $\Zp\ps{G}$-module $M$, denote by
$M[p^\infty]$ the $\Zp\ps{G}$-submodule of $M$ consisting of elements
of $M$ which are annihilated by some power of $p$. Howson \cite[Proposition
1.11]{Ho2}, and independently Venjakob \cite[Theorem 3.40]{V02}, showed that there is a
$\Zp\ps{G}$-homomorphism
\[ \varphi: M[p^\infty]\lra \bigoplus_{i=1}^s\Zp\ps{G}/p^{\alpha_i},\] whose
kernel and cokernel are pseudo-null $\Zp\ps{G}$-modules, and where
the integers $s$ and $\alpha_i$ are uniquely determined. The $\mu_G$-invariant of $M$ is then defined to be $\mu_G(M) = \displaystyle
\sum_{i=1}^s\alpha_i$.

We now extend the notion of torsion modules and pseudo-null modules to the case when $G$ is a general compact $p$-adic Lie group. A well-known theorem of Lazard asserts that $G$ contains an open normal subgroup $G_0$ which is pro-$p$ with no $p$-torsion (cf.\ \cite[Theorem 8.32]{DSMS}). By \cite[Proposition 5.4.17]{NSW}, we have
\[\Ext^i_{\Zp\ps{G}}(M,\Zp\ps{G}) \cong \Ext^i_{\Zp\ps{G_0}}(M,\Zp\ps{G_0})\]
for every finitely generated $\Zp\ps{G}$-module $M$. In view of this, we shall say that $M$ is a torsion $\Zp\ps{G}$-module (resp., psuedo-null $\Zp\ps{G}$-module) if $\Hom_{\Zp\ps{G}}(M,\Zp\ps{G})=0$ (resp., $\Ext^i_{\Zp\ps{G}}(M,\Zp\ps{G})=0$ for $i=0,1$). Equivalently, this is saying that $M$ is a torsion $\Zp\ps{G}$-module (resp., pseudo-null $\Zp\ps{G}$-module), whenever $M$ is a torsion $\Zp\ps{G_0}$-module (resp., pseudo-null $\Zp\ps{G_0}$-module) as in the preceding paragraph. (Also, compare with \cite[Discussion after Definition 2.6]{V02}).

\section{$K$-groups and Galois cohomology} \label{Arithmetic preliminaries}

We now describe a relation between the algebraic $K$-groups and Galois cohomology. To begin with,
$F$ will always denote a global function field of characteristic $\ell$. We shall, once and for all, fix a finite nonempty set $S$ of primes of $F$. Let $F_S$ be the maximal algebraic and separable extension of $F$ unramified outside $S$. For every extension $L$ of $F$ contained in $F_S$, write $G_S(L)$ for the Galois group $\Gal(F_S/L)$. Finally, we denote by $S(L)$ the set of primes of $L$ above $S$, and for brevity, we shall write $\Op_{L,S}$ for the ring of $S(L)$-integers of $L$.

For a ring $R$ with identity, write $K_n(R)$ for the algebraic $K$-groups of $R$ in the sense of Quillen \cite{Qui73a} (also see \cite{Kol, WeiKbook}). In this paper, we are interested in the group $K_n(\Op_{F,S})$. A theorem of Harder \cite[Korollar 3.2.3]{Ha} asserts that the group $K_n(\Op_{F,S})$ is finite for every $n\geq 2$. On the other hand, the group $K_n(F)$ is $p$-divisible by a result of Geisser-Levine \cite[Theorem 8.1]{GL}. By appealing to this and the localization sequence, one can then show that the group $K_n(\Op_{F,S})$ has order coprime to $\ell$ (we refer readers to \cite[Chap VI, Theorem 6.1 and Corollary 6.2]{WeiKbook} for details). Note that when $n=2$, this was earlier established by Bass-Tate \cite[Chapter II, Theorem 2.1]{BT}).

For the remainder of the paper, $p$ will always denote a prime $\neq \ell$. Denoting by $\mu_{p^n}$ the cyclic group generated by a primitive $p^n$-root of unity, we write $\mu_{p^\infty}$ for the direct limit of the groups $\mu_{p^n}$. The natural action of $G_S(F)$ on $\mu_{p^\infty}$ induces the cyclotomic character
\[\chi: G_S(F) \lra \mathrm{Aut}(\mu_{p^\infty}) \cong \Zp^{\times}.\]
For a discrete or compact $G_S(F)$-module $X$, we denote by $X(i)$ the $i$-fold Tate twist of $X$. More precisely, $X(i)$ is the $G_S(F)$-module which is $X$ as a $\Zp$-module but with a $G_S(F)$-action given by
\[ \sigma\cdot x = \chi(\sigma)^i\sigma x,  \]
where the action on the right is the original action of $G_S(F)$ on $X$.
Plainly, we have $X(0)=X$ and $\mu_{p^{\infty}} \cong \Qp/\Zp(1)$.
In \cite{Sou}, Soul\'e connected the higher $K$-groups with (continuous) Galois cohomology groups via the $p$-adic Chern class maps
\[ \mathrm{ch}_{i,k}^{(p)}: K_{2i-k}(\Op_{F,S})\ot \Zp \lra H^k(G_{S}(F), \Zp(i))\]
for $i\geq 2$ and $k =1,2$. (For the precise definition of these maps, we refer readers to loc.\ cit.) Soul\'e also proved that these maps are surjective (see \cite[Th\'eor\`{e}me 6(iii)]{Sou}; also see the work of Dwyer-Friedlander \cite{DF}). It is folklore (for instance, see \cite[Theorem 2.7]{Kol}) that the bijectivity of the $p$-adic Chern class maps follows as a consequence of the norm residue isomorphism theorem (previously also known as the Bloch-Kato(-Milnor) conjecture; see \cite{BK, Mil}), and the latter was proven by Rost and Voevodsky \cite{Vo} with the aid of a patch from Weibel \cite{Wei09}. As a consequence of these developments, we have the following identification between the Sylow $p$-subgroup of $K$-groups and Galois cohomology.

\bp \label{K2 = H2}
For $i\geq 2$, one has the following isomorphisms of finite abelian groups:
\[ K_{2i-1}(\Op_{F,S})[p^\infty] \cong H^1\big(G_{S}(F), \Zp(i)\big),\]
\[ K_{2i-2}(\Op_{F,S})[p^\infty] \cong H^2\big(G_{S}(F), \Zp(i)\big).\]
\ep

\section{Iwasawa cohomology groups} \label{Iwasawa cohomology groups}

We now introduce the Iwasawa cohomology groups which will be an important object for our study of $K$-groups. As before, $F$ will denote a global function field of characteristic $\neq p$, and we let $i\geq 2$. For every extension $\mathcal{L}$ of $F$ contained in $F_S$, the Iwasawa cohomology groups are defined by
 \[H^k_{\Iw,S}\big(\mathcal{L}/F, \Zp(i)\big):= \plim_L H^k\big(G_S(L), \Zp(i)\big), \]
 where $L$ runs through all finite extensions of $F$ contained in $\mathcal{L}$ and the transition maps are given by the corestriction maps. In the event that $\mathcal{L}/F$ is a finite extension, we plainly have
 $H^k_{\Iw,S}\big(\mathcal{L}/F, \Zp(i)\big)= H^k\big(G_S(\mathcal{L}), \Zp(i)\big)$.

Now suppose that $F_\infty$ is a $p$-adic Lie extension of $F$ contained in $F_S$. In other words, $F_\infty$ is a Galois extension of $F$, whose Galois group $G=\Gal(F_\infty/F)$ is a compact $p$-adic Lie group. We can now state and prove the first theorem of the paper.

\bt \label{H2 mu=0}
Let $i\geq 2$. Suppose that $F_\infty$ is a $p$-adic Lie extension of $F$ which is contained in $F_S$.  Then the following statements are valid.
\begin{enumerate}
  \item[$(a)$] The module $H^2_{\Iw,S}(F_\infty/F, \Zp(i))$ is torsion over $\Zp\ps{G}$. In the event that $F_\infty$ contains $F_\cyc$, we have $H^2_{\Iw,S}(F_\infty/F, \Zp(i))$ being finitely generated over $\Zp\ps{H}$, where $H= \Gal(F_\infty/F_\cyc)$.
  \item[$(b)$] $H^1_{\Iw,S}(F_\infty/F, \Zp(i))$ is a torsion $\Zp\ps{G}$-module. Furthermore, in the event that the group $G$ has dimension $\geq 2$, one has
     \[ H^1_{\Iw,S}(F_\infty/F, \Zp(i)) =0. \]
\end{enumerate}
\et

The remainder of the section is devoted to the proof of this theorem. As a start, we record the following version of Tate's descent spectral sequence for Iwasawa cohomology.

\bp \label{gal descent}
Suppose that $i\geq 2$.
Let $U$ be a closed normal subgroup of $G=\Gal(F_\infty/F)$ and write $\mathcal{L}$ for the fixed field of $U$. Then we have a homological spectral sequence
 $$ H_r\big(U, H^{-s}_{\Iw,S}(F_\infty/F, \Zp(i))\big)\Longrightarrow H^{-r-s}_{\Iw,S}\big(L_\infty/F, \Zp(i)\big) $$
and an isomorphism
\[ H^2_{\Iw,S}\big(F_\infty/F, \Zp(i)\big)_U \cong H^2_{\Iw,S}\big(\mathcal{L}/F, \Zp(i)\big).\]
Furthermore, if $U$ is an open subgroup of $G$, then $\mathcal{L}$ is a finite extension of $F$ and we have
\[ H^2_{\Iw,S}\big(F_\infty/F, \Zp(i)\big)_U \cong H^2\big(G_S(\mathcal{L}), \Zp(i)\big).\]
\ep

\bpf
When the group $G$ is commutative, the spectral sequence is a result of Nekov\'a\v{r} \cite[Proposition 4.2.3]{Ne}. For a noncommutative $G$, this follows from \cite[Proposition 1.6.5]{FK} or \cite[Theorem 3.1.8]{LimSh}. The isomorphisms of the proposition follows from reading off the initial $(0,-2)$-term of the spectral sequence.
\epf

We record another observation pertaining to the $p$-cohomological dimension of the group $G_S(F_\cyc)$ which will be an important component to our subsequent discussion.

\bl\label{cdpcyc} The $p$-cohomological dimension of the group $G_S(F_\cyc)$ is equal to $1$.
\el

\bpf
 Let $k$ be the constant field of $F$, and write $\bar{k}$ for the algebraic closure of $k$. Since the set $S$ is nonempty and $\mathrm{char}(F)\neq p$, it follows from \cite[Theorem 10.1.2(i)(c)]{NSW} that the group $G_S(\bar{k}F)$ is of $p$-cohomological dimension $1$. Since $G_S(\bar{k}F)$ is a closed subgroup of $G_S(F_{\mathrm{cyc}})$ with quotient $\Gal(\bar{k}F/F_{\mathrm{cyc}})$ whose profinite order is coprime to $p$, the two groups $G_S(F_{\mathrm{cyc}})$ and $G_S(\bar{k}F)$ have the same $p$-cohomological dimension (see \cite[Proposition 3.3.5(i)]{NSW}). The conclusion of the lemma now follows.
\epf

We are in position to prove Theorem \ref{H2 mu=0}.

\bpf[Proof of Theorem \ref{H2 mu=0}]
We first prove statement (a). The torsionness can be proved via a similar argument to that in \cite[Proposition 4.1.1]{LimKgrowth}. We will focus on the finite generation assertion.
 By base changing $F$, we may assume that $G$ is pro-$p$ with no $p$-torsion. Since
\[ H^2_{\Iw,S}\big(F_\infty/F, \Zp(i)\big)_H \cong H^2_{\Iw,S}\big(F_\cyc/F, \Zp(i)\big)\]
by Proposition \ref{gal descent}, a Nakayama's lemma argument reduces us to showing that $H^2_{\Iw,S}\big(F_\cyc/F, \Zp(i)\big)$ is finitely generated over $\Zp$. But since we also have
\[ H^2_{\Iw,S}\big(F_\cyc/F, \Zp(i)\big)/p \cong H^2_{\Iw,S}\big(F_\cyc/F, \Z/p(i)\big), \]
 this in turn reduces us to showing that $H^2_{\Iw,S}\big(F_\cyc/F, \Z/p(i)\big)$ is finite, or equivalently, $H^2_{\Iw,S}\big(F_\cyc/F, \Z/p(i)\big)$ is torsion over $\mathbb{F}_p\ps{\Gamma}$, where $\Ga = \Gal(F_\cyc/F)$. By considering the initial term of the spectral sequence
\[ \Ext^r_{\mathbb{F}_p\ps{\Ga}}\Big(H^{3-s}_{\Iw,S}\big(F_\cyc/F, \Z/p(i)\big),\mathbb{F}_p\ps{\Ga}\Big)\Longrightarrow  H^{3-r-s}\big(G_{S}(F_\cyc),\Z/p(i)\big)^{\vee}\]
(see \cite[1.6.12]{FK},  \cite[Theorem 4.5.1]{LimSh} or \cite[Theorem 5.4.5]{Ne}), we obtain an isomorphism
\[ \Ext^0_{\mathbb{F}_p\ps{\Ga}}\Big(H^{2}_{\Iw,S}\big(F_\cyc/F, \Z/p(i)\big),\mathbb{F}_p\ps{\Ga}\Big)\cong H^{2}\big(G_{S}(F_\cyc),\Z/p(i)\big)^{\vee} \]
where we note that the latter vanishes by virtue of Lemma \ref{cdpcyc}. Hence it follows that  $H^2_{\Iw,S}\big(F_\cyc/F, \Zp(i)\big)$ is torsion over $\mathbb{F}_p\ps{\Ga}$ as required. The proof of statement (a) is therefore complete.

We now turn to the verification of statement (b). As before, we may assume that $G$ is pro-$p$ with no $p$-torsion. Since the first assertion of statement (a) is telling us that
\[ \rank_{\Zp\ps{G}}\Big( H^2_{\Iw,S}\big(F_\infty/F, \Zp(i)\big) \Big)=0,\]
we therefore have
\[ -\rank_{\Zp\ps{G}} \Big(H^1_{\Iw,S}\big(F_\infty/F, \Zp(i)\big)\Big) = \sum_{s\geq 0}(-1)^s \rank_{\Zp\ps{G}}\Big(H^s_{\Iw,S}\big(F_\infty/F, \Zp(i)\big)\Big).\]
By a formula of Howson \cite[Theorem 1.1]{Ho}, the latter sum can be expressed as
\[ \sum_{r,s\geq 0}(-1)^{r+s} \rank_{\Zp}\Big(H_r\big(G, H^s_{\Iw,S}\big(F_\infty/F, \Zp(i)\big)\big)\Big). \]
Taking the spectral sequence of Proposition \ref{gal descent} into account, this sum is equal to
\[ \sum_{m\geq 0}(-1)^{m} \rank_{\Zp}\Big(H^m\big(G_S(F), \Zp(i)\big)\Big) =0, \]
where the final zero is a consequence of the cohomology groups appearing in the sum being finite. We have therefore established the first assertion of statement (b).

To continue, we shall work under the assumption that $G$ has dimension $\geq 2$. We now call upon another spectral sequence of Nekov\'a\v{r}
\[ \Ext^r_{\Zp\ps{G}}\Big(H^{s}\big(G_S(F_\infty), \Qp/\Zp(i)\big)^\vee,\Zp\ps{G}\Big)\Longrightarrow  H^{r+s}_{\Iw,S}\big(F_\infty/F,\Zp(i)\big)\]
(cf.\ \cite[1.6.12]{FK},  \cite[Theorem 4.5.1]{LimSh} or \cite[Theorem 5.4.5]{Ne}), whose low degree terms fit into the following exact sequence
\[ 0 \lra \Ext^1_{\Zp\ps{G}}\Big(H^{0}\big(G_S(F_\infty), \Qp/\Zp(i)\big)^\vee,\Zp\ps{G}\Big)\lra  H^{1}_{\Iw,S}\big(F_\infty/F,\Zp(i)\big) \]
\[  \lra \Ext^0_{\Zp\ps{G}}\Big(H^{1}\big(G_S(F_\infty), \Qp/\Zp(i)\big)^\vee,\Zp\ps{G}\Big).\]
As the group $G$ is of dimension $\geq 2$, it follows from a classical observation of Jannsen (cf.\ \cite[Corollary 2.6(b)(c)]{Jannsen89}) that $\Ext^1_{\Zp\ps{G}}\Big(M,\Zp\ps{G}\Big)=0$ whenever $M$ is finitely generated over $\Zp$. Since this is the case for the module $H^{0}\big(G_S(F_\infty), \Qp/\Zp(i)\big)^\vee$, we see that the leftmost term of the above exact sequence is zero.
Consequently, the module $H^{1}_{\Iw,S}\big(F_\infty/F,\Zp(i)\big)$, which is torsion over $\Zp\ps{G}$ as seen in the previous paragraph, injects into $\Ext^0_{\Zp\ps{G}}\Big(H^{1}\big(G_S(F_\infty), \Qp/\Zp(i)\big)^\vee,\Zp\ps{G}\Big)$ which is a reflexive $\Zp\ps{G}$-module (cf. \cite[Proposition 3.11(i)]{V02}). Thus, we must have $H^{1}_{\Iw,S}\big(F_\infty/F,\Zp(i)\big)=0$. This therefore completes the proof of the theorem.
\epf

\br
Some part of the assertions in Theorem \ref{H2 mu=0} should be generalized to more general Galois representations. Indeed, the theorem is inspired by a work of the author,  where, together with his coauthors, the $\Zp\ps{H}$-finite generation of the second Iwasawa cohomology group and vanishing of the first Iwasawa cohomology group for the Tate module of an abelian variety are established (see \cite{DKLL}). For a high-brow approach on the second assertion of Theorem \ref{H2 mu=0}(a), we refer readers to the work of Witte \cite{Wi}.
\er

We end with a byproduct of the discussion in this section. For the remainder of the section, we assume that $G$ is a compact $p$-adic Lie group with no $p$-torsion. The $G$-Euler characteristics of a $\Zp\ps{G}$-module $M$ is said to be well-defined if $H_j(G,M)$ is finite for every $j\geq 0$. In the event of such, the $G$-Euler characteristics of $M$ is then defined to be
\[ \chi(G,M) = \prod_{j\geq 0} \big(\# H_j(G,M)\big)^{(-1)^j}.\]
Note that since the group $G$ is assumed to have no $p$-torsion, it has finite $p$-cohomological dimension (see \cite[Corollaire 1]{Se65}), and therefore, the alternating product is a finite product. We now evaluate the $G$-Euler characteristics of the second Iwasawa cohomology group for $p$-adic Lie extension with dimension at least $2$.

\bt \label{Eulerchar}
Notation as before. Suppose further that $G$ has dimension $\geq 2$ and has no $p$-torsion. Then one has
\[ \chi\big(G,H^2_{\Iw}\big(F_\infty/F, \Zp(i)\big)\big) = \frac{\#K_{2i-2}(\Op_{F,S})[p^\infty]}{\#K_{2i-1}(\Op_{F,S})[p^\infty]}.\]
\et

\bpf
Indeed, from Theorem \ref{H2 mu=0}(b), we have $H^1_{\Iw}\big(F_\infty/F, \Zp(i)\big)=0$. Taking this latter observation into account, the spectral sequence
\[ H_r\big(G, H^{-s}_{\Iw}\big(F_\infty/F, \Zp(i)\big)\big)\Longrightarrow H^{-r-s}\big(G_S(F),\Zp(i)\big) \]
degenerates yielding
\[H_0\big(G,H^2_{\Iw}\big(F_\infty/F, \Zp(i)\big)\big) \cong K_{2i-2}(\Op_{F,S})[p^\infty],  \]
\[H_1\big(G,H^2_{\Iw}\big(F_\infty/F, \Zp(i)\big)\big) \cong K_{2i-1}(\Op_{F,S})[p^\infty]  \]
and $H_j\big(G,H^2_{\Iw}\big(F_\infty/F, \Zp(i)\big)\big) =0$ for $j\geq 2$. The conclusion is now immediate.
\epf

\br
Note that Proposition \ref{Eulerchar} does not require $F_\infty$ to contain $F_\cyc$.
\er

\section{Growth of even $K$-groups} \label{Growth results section}
In this section, we present growth formulas for even $K$-groups over certain $p$-adic Lie extensions of global function fields. For a positive integer $m$, we write $\ord_p(m)$ for the highest power of $p$ dividing $m$, i.e., $m = p^{\ord_{p}(m)}m'$ for some $p\nmid m'$. If $N$ is a finite group, we then write $\ord_p(N)$ for $\ord_p(|N|)$. In particular, if $N$ is a finite $p$-group, one has $|N| = p^{\ord_p(N)}$.

We begin with the cyclotomic $\Zp$-extension.

\bt \label{Zp asymptotic arith}
Let $i\geq 2$ be given. Denote by $F_n$ the intermediate subfields of $F_\cyc/F$ with $|F_n:F|=p^n$.
 Then we have
 \[ \ord_p\Big(K_{2i-2}(\Op_{F_n,S})[p^\infty]\Big) = \rank_{\Zp}\Big(H^{2}_{\Iw,S}\big(F_\infty/F,\Zp(i)\big)\Big) n + O(1).\]
\et

\bpf
Write $\Ga_n=\Gal(F_\infty/F_n)$. By a combination of Propositions \ref{K2 = H2} and \ref{gal descent}, we have
\[H^2_{\Iw,S}\big(F_\infty/F, \Zp(i)\big)_{\Gamma_n} \cong H^2\big(G_{S}(F_n), \Zp(i)\big)\cong K_{2i-2}(\Op_{F_n,S})[p^\infty].\]
The conclusion of the theorem now follows from an application of \cite[Proposition 5.3.17]{NSW} and taking Theorem \ref{H2 mu=0}(a) into account.
\epf

\br
In the number field context, the above was first proved by Coates \cite[Theorem 9]{C72} for $K_2$. Subsequently, Ji-Qin \cite[Theorem 3.3(2)]{JQ} extends Coates's result to the higher even $K$-groups. In these works, they have shown that the growth has the form of $\mu p^n+\la n + O(1)$. It is generally conjectured the leading term $\mu p^n$ vanishes for the cyclotomic $\Zp$-extension of a number field (see \cite{Iw73}; also see \cite[Subsection 5]{LimKgrowth}), although this vanishing is only known in the case when the base field is abelian over $\Q$ (see \cite{FW}). In the function field context, this vanishing is guaranteed by virtue of Theorem \ref{H2 mu=0}.
\er

Recall that Lazard's theorem (cf.\ \cite[Corollary 8.34]{DSMS}) tells us that every compact $p$-adic Lie group contains an open normal subgroup which is a uniform pro-$p$ group in the sense of \cite[Definition 4.1]{DSMS}. Therefore, one can always reduce consideration for a general compact $p$-adic Lie group to the case of a uniform pro-$p$ group which we will do for the remainder of this section. For a uniform pro-$p$ group, we shall write $G_n$ for the lower
$p$-series $P_{n+1}(G)$ which is defined recursively by $P_{1}(G) = G$, and
\[ P_{n+1}(G) = \overline{P_{n}(G)^{p}[P_{n}(G),G]}, ~\mbox{for}~ n\geq 1. \]
For every $n$, we have $|G:G_n| = |G:P_{n+1}(G)| = p^{dn}$ (see \cite[Sections 3 and 4]{DSMS}), where $d$ is the dimension of the uniform pro-$p$ group $G$. We now record two results which we require.

\bp \label{ZprtimeZp}
Let $G$ be a uniform pro-$p$ group which contains a closed subgroup $H\cong \Zp^{d-1}$ such that $G/H\cong \Zp$. Suppose that $M$ is a $\Zp\ps{G}$-module which is finitely generated over $\Zp\ps{H}$ and that $M_{G_n}$ is finite for every $n$. Then we have
\[ \ord_p(M_{G_n}[p^\infty]) = \rank_{\Zp\ps{H}}(M) np^{(d-1)n} +  O(p^{(d-1)n}).\]
\ep

\bpf
This proposition is implicit in the discussion in \cite{Lei, LiangL}, and a formal presentation of which is given in \cite[Proposition 2.4.1]{LimKgrowth}.
\epf

\bp \label{G/H gen}
Let $G$ be a uniform pro-$p$ group which contains a closed subgroup $H$ such that $G/H\cong \Zp$. Suppose that $M$ is a $\Zp\ps{G}$-module which is finitely generated over $\Zp\ps{H}$. Then we have
\[ \ord_p(M_{G_n}/p^n) \leq \rank_{\Zp\ps{H}}(M) np^{(d-1)n} + \mu_H(M)p^{(d-1)n}+ O(np^{(d-2)n}).\]
\ep

\bpf
This result can be thought as a variant of the results of Perbet in \cite{Per}, and its idea of proof follows closely to that in Perbet (see \cite[Proposition 2.4]{LimUpper} for the details).
\epf

Building on our Theorem \ref{H2 mu=0}, we can prove the following asymptotic formula for the growth of even $K$-groups for a special class of $p$-adic Lie extension.

\bt \label{uniform asymptotic arith2}
Let $i\geq 2$ be given.
  Suppose that $F_\infty$ is a $p$-adic Lie extension of $F$ containing $F_\cyc$ with $\Gal(F_\infty/F_\cyc)\cong \Zp^{d-1}$. Let $F_n$ be the intermediate subfield of $F_\infty/F$ with $\Gal(F_\infty/F_n) = G_n$.
 Then we have
 \[ \ord_p\Big(K_{2i-2}(\Op_{F_n,S})[p^\infty]\Big) = \rank_{\Zp\ps{H}}\Big(H^2_{\Iw,S}\big(F_\infty/F, \Zp(i)\big)\Big) np^{(d-1)n} +  O(p^{(d-1)n}). \]
\et

\bpf
 By Theorem \ref{H2 mu=0}, the module $H^2_{\Iw,S}\big(F_\infty/F, \Zp(i)\big)$ is finitely generated over $\Zp\ps{\Gal(F_\infty/F_\cyc)}$. Therefore, we may combine Propositions \ref{gal descent} and \ref{ZprtimeZp} to obtain the asymptotic formula.
\epf

For a general $p$-adic Lie extension which contains $F_\cyc$, we at least have an upper bound.

\bp \label{uniform asymptotic arith3}
Let $i\geq 2$ be given.
 Suppose that $F_\infty$ is a $p$-adic Lie extension of $F$ containing $F_\cyc$ with $G$ being a uniform pro-$p$ group of dimension $d$. Let $F_n$ be the intermediate subfield of $F_\infty/F$ with $\Gal(F_\infty/F_n) = G_n$.
 Then we have
 \[ \ord_p\Big(K_{2i-2}(\Op_{F_n})[p^n]\Big) \leq \rank_{\Zp\ps{H}}\Big(H^2_{\Iw,S}\big(F_\infty/F, \Zp(i)\big)\Big) np^{(d-1)n} \hspace{1.5in}\]
 \[\hspace{1.5in} + ~\mu_{H}\Big(H^2_{\Iw,S}\big(F_\infty/F, \Zp(i)\big)\Big) p^{(d-1)n} +  O(np^{(d-2)n}). \]
 \ep

\bpf
 In this context, we make use of Proposition \ref{G/H gen} and note that $\#K_{2i-2}(\Op_{F_n})/p^n = \#K_{2i-2}(\Op_{F_n})[p^n]$ in view of the finiteness of $K_{2i-2}(\Op_{F_n})$.
\epf

\section{Comparing limits} \label{Klimitthm section}

We come to the final section of the paper, whose goal is to prove a relation between the direct limit and the inverse limit of even $K$-groups.  For two finite subextensions $L\subseteq L'$ of $F$ contained in $F_\infty$, the inclusion $\Op_{L,S}\lra \Op_{L',S}$ induces a map $\jmath_{L/L'}: K_{2i-2}(\Op_{L,S})[p^\infty]\lra K_{2i-2}(\Op_{L',S})[p^\infty]$ by functoriality. (Note that the induced map may not be injective in general.) In the other direction, there is the norm map (also called the transfer map) $\Tr_{L'/L}: K_{2i-2}(\Op_{L',S})[p^\infty]\lra K_{2i-2}(\Op_{L,S})[p^\infty]$. Consider the following direct limit and inverse limit
\[  \ilim_L K_{2i-2}(\Op_{L,S})[p^\infty] \quad \mbox{and}\quad \plim_L K_{2i-2}(\Op_{L,S})[p^\infty],\]
whose transition maps are given by the maps $\jmath_{L/L'}$ and $\Tr_{L'/L}$ respectively.
The main result of this section is as follows.

\bt \label{Klimitthm}
Let $i\geq 2$ and let $F_\infty$ be a $p$-adic Lie extension of $F$ which is contained in $F_S$. Then there is an isomorphism
\[  \left(\ilim_L K_{2i-2}(\Op_{L,S})[p^\infty]\right)^\vee \cong \Ext^1_{\Zp\ps{G}}\left(\plim_L K_{2i-2}(\Op_{L,S})[p^\infty], \Zp\ps{G}\right)\]
of $\Zp\ps{G}$-modules.
\et

The following is an immediate corollary of the theorem.

\bc
Retain settings as in Theorem \ref{Klimitthm}. Then we have $\ilim_L K_{2i-2}(\Op_L)[p^\infty]=0$ if and only if
$\plim_L K_{2i-2}(\Op_L)[p^\infty]$ is pseudo-null over $\Zp\ps{G}$.
\ec

We now give the proof of Theorem \ref{Klimitthm}.

\bpf[Proof of Theorem \ref{Klimitthm}] Now, if the dimension of $G$ is $1$, then by base-changing, we may assume that $F_\infty/F$ is a $\Zp$-extension, which is a commutative $p$-adic Lie extension. In this situation, the theorem follows from a direct similar argument to that in \cite[Theorem 3.1]{LimKlimit}. Therefore, we may assume that $G$ has dimension $\geq 2$. From the discussion in \cite[Chap. III]{Sou}, there is a commutative diagram
\begin{equation*} \label{N cor}\entrymodifiers={!! <0pt, .8ex>+} \SelectTips{eu}{}
\xymatrix{
      K_{2i-2}(\Op_{L,S})[p^\infty] \ar[r]^{\mathrm{ch}_i^L} \ar[d]^{\Tr_{L/F}} &  H^2\left(G_S(L), \Zp(i)\right) \ar[d]^{\mathrm{cor}} \\
      K_{2i-2}(\Op_{F,S})[p^\infty] \ar[r]^{\mathrm{ch}_i^F} &  H^2\left(G_S(F), \Zp(i)\right)}
      \end{equation*}
which in turn yields
\[ H^2_{\Iw}\big(F_\infty/F, \Zp(i)\big) \cong  \plim_{L} K_{2i-2}(\Op_{L,S})[p^\infty]. \]
For brevity, write $M=H^2_{\Iw}\big(F_\infty/F, \Zp(i)\big)$.
By Jannsen's result (cf.\ \cite[Theorem 2.1]{Jannsen89} or \cite[Theorem 5.4.13]{NSW}), there is a short exact sequence
\begin{equation}\label{Jannsenses}
 0 \lra \Big(\ilim_U(M_U[p^\infty])\Big)^\vee \lra \Ext^1_{\Zp\ps{G}}(M,\Zp\ps{G})
  \lra \Big(\ilim_U H_1(U,M)\ot\Qp/\Zp\Big)^\vee \lra 0.
\end{equation}
where the transition map in the limit is given by norm map given as follow: For two open normal subgroups $V\subseteq U$ of $G$, the norm map $N_{U/V}$ on $M_V$ factors through $M_U=(M_V)_{U/V}$ to yield a map $M_U\lra M_V$, which by abuse of notation is also denoted by $N_{U/V}$.

By a similar argument to that in \cite[Lemma 3.6]{LimKlimit}, we obtain the following commutative diagram
\[ \entrymodifiers={!! <0pt, .8ex>+} \SelectTips{eu}{}\xymatrix{
      H^2_{\Iw}\big(F_\infty/F, \Zp(i)\big)_U \ar[r]_{\sim}\ar[d]_{N_{U/V}} &  H^2\big(G_S(L), \Zp(i)\big) \ar[d]_{\mathrm{res}}&  K_{2i-2}(\Op_{L,S})[p^\infty] \ar[l]_(.45){\mathrm{ch}_i^{L}}^(.45){\sim} \ar[d]_{j_{L/L'}} \\
      H^2_{\Iw}\big(F_\infty/F, \Zp(i)\big)_V  \ar[r]_{\sim} &  H^2\big(G_S(L'), \Zp(i)\big) & \ar[l]_(.45){\mathrm{ch}_i^{L’}}^(.45){\sim} K_{2i-2}(\Op_{L',S})[p^\infty]
     }\]
which yields an identification of the leftmost term in the short exact sequence (\ref{Jannsenses}) with
\[  \ilim_{L} K_{2i-2}(\Op_{L,S})[p^\infty].\]
It therefore remains to show the rightmost term in the short exact sequence (\ref{Jannsenses}) vanishes. For this, it suffices to show that $H_1(U,M)$ is finite for every open subgroup $U$ of $G$. But as already seen in the proof of Theorem \ref{Eulerchar}, if we denote by $L$ the fixed field of $U$, then one has
\[H_1\Big(U,H^2_{\Iw}\big(F_\infty/F, \Zp(i)\big)\Big) \cong K_{2i-1}(\Op_{L,S})[p^\infty],  \]
where the latter $K$-group is finite by Harder's theorem. The proof of our theorem is now complete.
\epf

\end{document}